\documentstyle{article}

\newcommand{\nc}{\newcommand}
\nc{\be}{\begin{equation}}
\nc{\ee}{\end{equation}}
\nc{\dfrac}{\displaystyle \frac }
\begin{document}

\author{A. Yildiz \\
Feza G\"ursey Institute, P.O. Box 6, 81220, \c Cengelk\"oy, Istanbul, Turkey}
\title{Braided Oscillators }
\maketitle

\begin{abstract}
A generalized oscillator algebra is proposed and the braided Hopf algebra
structure for this generalized oscillator is investigated. Using the
solutions for the braided Hopf algebra structure, two types of braided
Fibonacci oscillators are introduced. This leads to two types of braided
Biedenharn-Macfarlane oscillators as special cases of the Fibonacci
oscillators. We also find the braided Hopf algebra solutions for the three
dimensional braided space. One of these, as a special case, gives the Hopf
algebra given in the literature.
\end{abstract}

\pagebreak

\section{Introduction }

 \baselineskip=27pt The harmonic oscillator has a wide variety of
applications from quantum optics to the realizations of the
angular momentum algebra and hence the deformations of the
oscillator algebra play
an important role in $q$-deformed theories. The realization of the $q$%
-deformed angular momentum algebra by Biedenharn-Macfarlane oscillators$^{1}$
and the realization of the Hermitian braided matrices by a pair of $q$%
-oscillators$^{2}$ are some of the examples. The two parameter deformations
and some of their applications can be found in$^{3}$. Braided group theory
(a self contained review can be found in$^{4}$) deforms the notion of tensor
product (called braided tensor product) and hence deforms the independence
of the objects. Although braided groups arise in the formulation of quantum
group covariant structures, the idea of braiding can be used without any
reference to quantum groups to generalize the statistics$^{5}$.

The permutation map $\pi$ ($\pi :A\otimes B\rightarrow B\otimes A$) in the
tensor product algebra of boson algebras

$(a\otimes b)(c\otimes d)=a\pi (b\otimes c)d=ac\otimes bd$

\noindent is replaced by a generalized map called braiding $\psi $ ($\psi
:A\otimes B\rightarrow B\otimes A$) such that

$(a\otimes b)(c\otimes d)=a\psi (b\otimes c)d.$

This generalization leads to the generalization of the Hopf algebra called
braided Hopf algebra$^6$ whose axioms in algebraic (not diagrammatic) form
read as

\begin{eqnarray}  \label{eq:bhopf}
m\circ (id\otimes m)&=&m\circ (m\otimes id)  \nonumber \\
m\circ (id\otimes \eta )&=&m\circ (\eta \otimes id) = id  \nonumber \\
(id\otimes \Delta )\circ \Delta &=&(\Delta \otimes id)\circ \Delta  \nonumber
\\
(\epsilon \otimes id)\circ \Delta &=&(id\otimes \epsilon )\circ \Delta = id
\nonumber \\
m\circ (id\otimes S)\circ \Delta &=&m\circ (S\otimes id)\circ \Delta =\eta
\circ \epsilon  \nonumber \\
\psi \circ (m \otimes id) &=& (id \otimes m)\circ (\psi \otimes id)\circ (id
\otimes \psi)  \nonumber \\
\psi \circ (id \otimes m) &=& (m \otimes id)\circ (id \otimes \psi)\circ
(\psi \otimes id)  \nonumber \\
(id \otimes \Delta)\circ \psi &=& (\psi \otimes id)\circ (id \otimes \psi
)\circ (\Delta \otimes id)  \nonumber \\
(\Delta \otimes id)\circ \psi &=& (id \otimes \psi )(\psi \otimes id)\circ
(id \otimes \Delta)  \nonumber \\
\Delta \circ m &=& (m \otimes m)(id\otimes \psi \otimes id)\circ (\Delta
\otimes \Delta ) \\
S\circ m &=& m\circ \psi \circ (S\otimes S)  \nonumber \\
\Delta \circ S &=& (S\otimes S)\circ \psi \circ \Delta  \nonumber \\
\epsilon \circ m &=& \epsilon \otimes \epsilon  \nonumber \\
(\psi \otimes id)\circ (id \otimes \psi)\circ (\psi \otimes id)&=&(id
\otimes \psi)\circ (\psi \otimes id)\circ (id \otimes \psi)  \nonumber
\end{eqnarray}

\noindent where $m:A\otimes A\rightarrow A$ is the multiplication map, $%
\Delta :A\rightarrow A\otimes A$ is the comultiplication map, $\eta
:K\rightarrow A$ is the unit map, $\epsilon :A\rightarrow K$ is the counit
map, $S:A\rightarrow A$ is the antipode map and $\psi :A\otimes A\rightarrow
A\otimes A$ is the braiding map. The consistency of the relations (\ref
{eq:bhopf}) requires that

\be
\Delta (1_A)=1_A\otimes 1_A,\ \ \psi (1_A\otimes a)=a\otimes 1_A,\ \ \psi
(a\otimes 1_A)=1_A\otimes a \ \ \forall a\in A \ee
\noindent where $1_A$ is the identity of the algebra A. From now on we will
drop the subscript and write $1$ for the identity of the algebra. All these
maps are linear. Note that in the limit $\psi \rightarrow \pi $ the braided
Hopf algebra axioms reduce to the ordinary Hopf algebra axioms. The braided
Hopf algebra axioms$^7$ reduce to the axioms given above when the counit map
($\epsilon$) is an algebra homomorphism. The $*$-structure for a braided
algebra $B$ is different$^8$ from the non-braided one such that

\begin{eqnarray}
\Delta \circ * &=& \pi \circ (* \otimes *)\circ \Delta  \nonumber \\
S\circ * &=& *\circ S \\
(a\otimes b)^{*} &=& b^{*}\otimes a^{*},\ \forall a,b\in B.  \nonumber
\end{eqnarray}

\section{Generalized Oscillator }
We propose a generalized oscillator algebra generated by $a,\ a^{*},\ q^N$
and $1$ satisfying

\begin{eqnarray}  \label{eq:go1}
aq^N &=& qq^Na,  \nonumber \\
q^Na^{*} &=& qa^{*}q^N, \\
aa^{*} -Q_1a^{*}a &=& Q_2q^{2N}+Q_3q^N+Q_4  \nonumber
\end{eqnarray}

\noindent where $q,Q_1,Q_2,Q_3,Q_4$ are real constants whose values
determine
the type of the oscillator. For example, if $Q_1$ is a free parameter then
the $Q_2=1,\ Q_3=Q_4=0$ case and the $Q_3=1,\ Q_2=Q_4=0$ case define two
different Fibonacci oscillators. For the $*$-structure we impose $%
(q^N)^{*}=q^N,\ \ (a^{*})^{*}=a.$ The actions of the generators on the
Hilbert space are given by

\begin{eqnarray}  \label{eq:actions1}
a\mid n\rangle & = & a_n\mid n-1\rangle  \nonumber \\
a^{*}\mid n\rangle & = & a_{n+1}^{*}\mid n+1\rangle \\
q^N\mid n\rangle & = & q^n\mid n\rangle.  \nonumber
\end{eqnarray}

Using the fact that for a given algebra the Hopf algebra structure is not
unique, we write the general forms of the coproducts
\begin{eqnarray}
\Delta (q^N) &=& A_1q^N\otimes q^N+A_2a\otimes a^{*} +A_3a^{*}\otimes
a+A_41\otimes q^N+A_5q^N\otimes 1 +A_61\otimes1,  \nonumber \\
\Delta (a) &=&B_1q^N\otimes a+B_2a\otimes q^N+B_31\otimes a+B_4a\otimes 1, \\
\Delta (a^{*}) &=& B_1a^{*}\otimes q^N+B_2q^N\otimes a^{*} +B_3a^{*}\otimes
1+B_41\otimes a^{*},  \nonumber
\end{eqnarray}

\noindent the counits

\be
\epsilon (q^N)=e_1, \ \ \ \epsilon (a)=\epsilon (a^{*})=e_2, \ee
\noindent the antipodes
\begin{eqnarray}
S (q^N) &=& k_1q^N+k_2a+k_3a^{*}+k_4,  \nonumber \\
S (a) &=& m_1q^N+m_2a+m_3a^{*}+m_4, \\
S (a^{*}) &=&m_1q^N+m_2a^{*}+m_3a+m_4  \nonumber
\end{eqnarray}

\noindent and the braidings

\begin{eqnarray}
\psi (q^N \otimes q^N) &=& g_1q^N\otimes q^N+ g_2a\otimes a^{*}
+g_3a^{*}\otimes a+g_41\otimes q^N+g_5q^N+g_61\otimes1,  \nonumber \\
\psi(q^N\otimes a) &=& d_1a\otimes q^N+d_2q^N\otimes a+d_31\otimes a
+d_4a\otimes 1,  \nonumber \\
\psi (a^{*}\otimes q^N) &=& d_1q^N\otimes a^{*}+ d_2a^{*}\otimes q^N
+d_3a^{*}\otimes 1+d_41\otimes a^{*},  \nonumber \\
\psi(q^N\otimes a^{*}) &=& f_1a^{*}\otimes q^N +f_2q^N\otimes
a^{*}+f_31\otimes a^{*}+f_4a^{*}\otimes 1,  \nonumber \\
\psi (a\otimes q^N) &=& f_1q^N\otimes a+f_2a\otimes q^N+f_3a\otimes 1+
f_41\otimes a,  \nonumber \\
\psi (a\otimes a) &=& z_1a\otimes a, \\
\psi (a^{*}\otimes a^{*}) &=& z_1a^{*}\otimes a^{*},  \nonumber \\
\psi(a\otimes a^{*}) &=& b_1q^N\otimes q^N+b_2a\otimes a^{*}
+b_3a^{*}\otimes a+b_41\otimes q^N+b_5q^N\otimes 1+b_61\otimes1,
\nonumber
\\
\psi(a^{*}\otimes a) &=& c_1q^N\otimes q^N+c_2a\otimes a^{*}
+c_3a^{*}\otimes a+c_41\otimes q^N+c_5q^N\otimes 1+c_61\otimes 1  \nonumber
\end{eqnarray}

\noindent where symbols with a subscript are the constants to be determined.

\section{Braided Hopf Algebra Structure of the Generalized Oscillator }

To find the general braided Hopf algebra structure for the
oscillator algebra given by (\ref{eq:go1}) we substitute these
general forms into the braided Hopf algebra axioms
(\ref{eq:bhopf}) and find the solutions using the computer
programming Mapple V. The constants which are the same for all
solutions read as

\be\label{eq:constants1} Q_4=A_1=A_6=B_1=B_2=0,\ A_4=A_5=B_3=B_4=1, \ee
\be \label{eq:constants2} k_2=k_3=k_4=m_1=m_3=m_4=e_1=e_2=0,\ m_2=-1, \ee
\be \label{eq:constants3}
b_4=b_5=b_6=c_4=c_5=c_6=g_4=g_5=g_6=d_3=d_4=f_3=f_4=0. \ee

\noindent The solutions for the other parameters are given in the tables.
The constants given by (\ref{eq:constants1})-(\ref{eq:constants3}) show that
the antipodes and the counits of all generators and the coproducts of
raising/lowering operators are uniquely determined. We also found that for a
free deformation parameter $Q_1$ at most one of the other deformation
parameters, namely $Q_2$ or $Q_3$, is nonzero. For $Q_2\neq 0$ , we have

\be \label{eq:q2} aa^{*}-Q_1a^{*}a=Q_2q^{2N}. \ee

\noindent Without loss of generality we can take $Q_2=1$ or by rescaling $a$
and $a^{*}$ the oscillator relation (\ref{eq:q2}) can be reduced to

\be \label{eq:q21} aa^{*}-Q_1a^{*}a=q^{2N}. \ee

\noindent We found that there are only six possible braided Hopf algebra
solutions for this two parameter oscillator which are given in Table \ref
{table1}. Defining the free deformation parameter $Q_1\equiv p^{-2}$ the
algebra (\ref{eq:q21}) can be rewritten as

\be \label{eq:q21p} aa^{*}-p^{-2}a^{*}a=q^{2N} \ee

\noindent which is the more familiar form of the $q,p$ oscillator algebra
(Fibonacci oscillator) which we call the first-type Fibonacci oscillator.
Using the actions (\ref{eq:actions1}) the eigenvalue of the operator $a^{*}a$
on the state $\mid n\rangle$ is found to be

\be
a_n^{*}a_n=\displaystyle \frac{p^{-2n}-q^{2n}}{p^{-2}-q^2}. \ee

\noindent Substituting $Q_1=q^{-2}$ into (\ref{eq:q21}) we obtain the
Biedenharn-Macfarlane oscillator algebra which we call the first-type
Biedenharn-Macfarlane oscillator defined by

\be\label{eq:bm1} aa^{*}-q^{-2}a^{*}a=q^{2N}. \ee

\noindent The eigenvalue of the operator $a^{*}a$ on the state $\mid
n\rangle $ is found to be

\be
a_n^{*}a_n=\displaystyle \frac{q^{2n}-q^{-2n}}{q^{2}-q^{-2}}. \ee

\noindent There are only six braided Hopf algebra solutions for
the first-type Biedenharn Macfarlane oscillator (\ref{eq:bm1})
which can be obtained by substituting $Q_1=q^{-2}$ into the
solutions given in Table \ref {table1}. \noindent Similarly, for
$Q_3\neq 0$ the oscillator relation

\be \label{eq:q3} aa^{*}-Q_1a^{*}a=Q_3q^{N} \ee

\noindent can again be reduced to

\be \label{eq:q31} aa^{*}-Q_1a^{*}a=q^{N} \ee

\noindent which we call the second-type Fibonacci oscillator. Setting the
free deformation $Q_1\equiv p^{-1}$ the eigenvalue of the operator $a^{*}a$
on the state $\mid n\rangle$ is found to be

\be
a_n^{*}a_n=\displaystyle \frac{p^{-n}-q^{n}}{p^{-1}-q}. \ee

The second-type Fibonacci oscillator (\ref{eq:q31}) has only two braided
Hopf algebra solutions which are given in Table \ref{table3}. The
second-type Biedenharn-Macfarlane oscillator can be obtained simply by
substituting $Q_1=q^{-1}$ into (\ref{eq:q31}) and the same substitution into
Table \ref{table3} gives braided Hopf algebra solutions.

A wide variety of one parameter oscillators can be obtained by assigning $%
Q_1=f(q)$ in the algebras (\ref{eq:q21}) and (\ref{eq:q31}). The braided
Hopf algebra solutions for these oscillators can be obtained by substituting
$Q_1=f(q)$ into the Table \ref{table1} and Table \ref{table3}. However,
there are extra solutions for some values $Q_1$ which are given in Table \ref
{table5}.

The braided Hopf algebra structure of the quantum space (called braided
space) defined by

\be\label{eq:qs} x_ix_j=qx_jx_i\ \ i>j \ee

\noindent determines the structure of the braided integration,
derivation and Fourier transform defined on that space$^{9}$. For
$Q_1=q,\ Q_2=0,\ Q_3=0,\ Q_4=0$ the algebra (\ref{eq:go1}) reduces
to the three dimensional quantum space with the identifications

\be
x_1\equiv a^{*},\ \ x_2\equiv q^{N},\ \ x_3\equiv a. \ee

\noindent The braided Hopf algebra solutions for the three dimensional
braided space are given in Table \ref{table4}. Setting the free parameter $%
g_1=q^2$ in the sol5 of Table \ref{table4} gives the braided Hopf
algebra given in the literature$^{10}$ and references therein) as
a special case.

\section{Conclusion}
The braidings imply the relations between independent copies of
algebras. For example, the implication of the braiding $\psi
(a\otimes a)$ can be found by using the identifications $a_1\equiv
a\otimes 1$ and $a_2\equiv 1\otimes a$ such that

$a_2a_1=(1\otimes a)(a\otimes 1)=\psi (a\otimes a)=za\otimes
a=za_1a_2$.

Thus the braiding relations imply a 2-body system of oscillators
which can be extended to $n$-body oscillators using the $n$-fold
braided tensor product as done by Baskerville and Majid in the
context of the braided version of the $q$-Heisenberg
algebra$^{11}$. The Fock space representations of the $n$-fold
braided tensor product of oscillators can also be found. This
requires the braidings of the generators of the algebra with the
states. All these constructions depend on the braiding relations
of the generators of the algebra. Since each solution for the
braiding gives a different system of oscillators, the solutions we
found may provide a general framework for the interacting
oscillators and hence for the statistical mechanical quantities
calculated by using these oscillators.

Because of the connection between symmetry and statistics it is
interesting to investigate the underlying symmetry transformations of the
braided oscillators which gives rise to the braiding relations
given in the tables. It may also be interesting to find the unbraiding
transformations$^{12}$ or the decoupling of the braided
oscillators. The generalization of supersymmetry to fractional
supersymmetry requires the deformation parameter to be a root of
unity$^{13}$ and this case deserves to be discussed on a seperate study.

The braided Hopf algebra solutions we present for the oscillators and for
the three dimensional braided space may provide a general frame on which
other structures can be defined.

 \pagebreak

{\bf References}

\begin{description}
\item  {$^1$}L. C. Biedenharn, J. Phys. A: Math. Gen. {\bf 22}, L873 (1989),
A. J. Macfarlane, J. Phys. A: Math. Gen. {\bf 22}, 4581 (1989)

\item  {$^2$}M. Arik  and A. Yildiz,   J. Phys. A: Math. Gen. {\bf 30},
L255 (1997)

\item  {$^3$}M. Arik, E. Demircan, T. Turgut, L. Ekinci and M. Mungan,
Z. Phys. C {\bf 55}, 89 (1992), R. Chakrabarti  and R. Jagannathan
  J. Phys. A: Math. Gen. {\bf 24}, L711 (1991), M. Daoud  and
M. Kibler,  Phys. Lett. A {\bf 206}, 13 (1995)

\item  {$^4$}S. Majid,  Beyond supersymmetry and quantum symmetry (an
introduction to braided groups and braided matrices) {\it Quantum
Groups, Integrable Statistical Models and Knot Theory} edited by
M. L. Ge and H. J. Vega (Singapore: World Scientific), pp.231-282
(1993)

\item  {$^5$}R. S. Dunne, A.J. Macfarlane, J.A. Azcarraga and J.C.
Bueno,  Phys. Lett. B {\bf 387}, 294 (1996)

\item  {$^6$}S. Majid,   Algebras and Hopf algebras in braided categories
{\it Lecture Notes in Pure and Applied Mathematics} vol 158 (New
York: Marcel Dekker), pp. 55-105 (1994)

\item  {$^7$}M. Durdevich,   Israel Journal of Mathematics {\bf 98},
329 (1997)

\item  {$^8$}Majid S  {\it Foundations of Quantum Group Theory}
(Cambridge: Cambridge University Press), (1995)

\item  {$^{9}$}A. Kempf  and S. Majid,  J. Math. Phys. {\bf 35},
6802 (1994)

\item  {$^{10}$}G. Carnovale,  J. Math. Phys. {\bf 40},
5972 (1999)

\item  {$^{11}$}W. K. Baskerville  and S. Majid,    J. Math. Phys. {\bf 34},
3588 (1993)

\item {$^{12}$}G. Fiore, H. Steinacker and J. Wess, Mod. Phys.
Lett. A {\bf 16}, 261 (2001)

\item {$^{13}$}R.S. Dunne, J. Math. Phys. {\bf 40}, 1180 (1999),
H. Ahmedov and O. F. Dayi, Mod. Phys. Lett. A {\bf 15}, 1801 (2000),
H. Ahmedov, A. Yildiz and Y. Ucan, J. Phys. A {\bf 34}, 6413 (2001)
\end{description}

 \pagebreak
\begin{table}[tbp]
\caption{Braided Fibonacci oscillator of the first type}
\label{table1}{\scriptsize
\begin{tabular}{|c||c|c|c|c|c|c|}
\hline
\multicolumn{7}{|c|}{$aa^{*}-Q_1a^{*}a=q^{2N}$,\ $aq^{N}=qq^Na$,\ $%
q^{N}a^{*}=qa^{*}q^{N}$} \\ \hline
& sol1 & sol2 & sol3 & sol4 & sol5 & sol6 \\ \hline
$k_1=$ & $-1$ & $-1$ & $-1$ & $-1$ & $-1$ & $-1$ \\ \hline
$b_1=$ & $1$ & $\displaystyle 1-\frac{q}{\sqrt{Q_1}}$ & $1+\displaystyle
\frac{q}{\sqrt{Q_1}}$ & $\displaystyle \frac{Q_1}{q^2}$ & $0$ & $0$ \\ \hline
$b_2=$ & $\displaystyle \frac{q^2-Q_1}{Q_1}$ & $\displaystyle \frac{%
(q^2-Q_1)(\sqrt{Q_1}+q)}{qQ_1}$ & $\displaystyle \frac{(Q_1-q^2)(\sqrt{Q_1}%
-q)}{qQ_1}$ & $0$ & $0$ & $0$ \\ \hline
$b_3=$ & $Q_1$ & $Q_1$ & $Q_1$ & $\displaystyle \frac{Q_1^2}{q^2}$ & $Q_1$
&
$Q_1$ \\ \hline
$z=$ & $\displaystyle \frac{q^2}{Q_1}$ & $\displaystyle \frac{q^2}{Q_1}$ & $%
\displaystyle \frac{q^2}{Q_1}$ & $\displaystyle \frac{Q_1}{q^2}$ & $%
\displaystyle \frac{Q_1}{q^2}$ & $\displaystyle \frac{Q_1}{q^2}$ \\ \hline
$c_1=$ & $-\displaystyle \frac{q^2}{Q_1^2}$ & $0$ & $0$ & $\displaystyle
\frac{-1}{Q_1}$ & $-\displaystyle \frac{q+\sqrt{Q_1}}{qQ_1}$ & $%
\displaystyle \frac{\sqrt{Q_1}-q}{qQ_1}$ \\ \hline
$c_2=$ & $\displaystyle \frac{q^2}{Q_1^2}$ & $\displaystyle \frac{1}{Q_1}$ &
$\displaystyle \frac{1}{Q_1}$ & $\displaystyle \frac{1}{Q_1}$ & $%
\displaystyle \frac{1}{Q_1}$ & $\displaystyle \frac{1}{Q_1}$ \\ \hline
$c_3=$ & $0$ & $0$ & $0$ & $\displaystyle \frac{Q_1-q^2}{q^2}$ & $%
\displaystyle \frac{(q^2-Q_1)(q-\sqrt{Q_1})}{q^2\sqrt{Q_1}}$ & $%
\displaystyle \frac{(Q_1-q^2)(q+\sqrt{Q_1})}{q^2\sqrt{Q_1}}$ \\ \hline
$d_1=$ & $\displaystyle \frac{q}{Q_1}$ & $\displaystyle \frac{q}{Q_1}$ & $%
\displaystyle \frac{q}{Q_1}$ & $q^{-1}$ & $q^{-1}$ & $q^{-1}$ \\ \hline
$d_2=$ & $0$ & $0$ & $0$ & $\displaystyle \frac{Q_1-q^2}{q^2}$ & $%
\displaystyle \frac{Q_1-q^2}{q^2}$ & $\displaystyle \frac{Q_1-q^2}{q^2}$ \\
\hline
$f_1=$ & $q$ & $q$ & $q$ & $\displaystyle \frac{Q_1}{q}$ & $\displaystyle
\frac{Q_1}{q}$ & $\displaystyle \frac{Q_1}{q}$ \\ \hline
$f_2=$ & $\displaystyle \frac{q^2-Q_1}{Q_1}$ & $\displaystyle \frac{q^2-Q_1}{%
Q_1}$ & $\displaystyle \frac{q^2-Q_1}{Q_1}$ & $0$ & $0$ & $0$ \\ \hline
$g_1=$ & $\displaystyle \frac{q^2}{Q_1}$ & $-\displaystyle \frac{q}{\sqrt{Q_1%
}}$ & $\displaystyle \frac{q}{\sqrt{Q_1}}$ & $\displaystyle \frac{Q_1}{q^2}$
& $\displaystyle \frac{\sqrt{Q_1}}{q}$ & $\displaystyle - \frac{\sqrt{Q_1}}{q%
}$ \\ \hline
$g_2=$ & $0$ & $\displaystyle \frac{(q^2-Q_1)(\sqrt{Q_1}+q)}{qQ_1}$ & $%
\displaystyle \frac{(Q_1-q^2)(\sqrt{Q_1}-q)}{qQ_1}$ & $0$ & $0$ & $0$ \\
\hline
$g_3=$ & $0$ & $0$ & $0$ & $0$ & $-\displaystyle \frac{Q_1(Q_1-q^2)^2}{%
q^2(Q_1+q\sqrt{Q_1})}$ & $\displaystyle \frac{Q_1(Q_1-q^2)^2}{q^2(q\sqrt{Q_1}%
-Q_1)}$ \\ \hline
\multicolumn{7}{|c|}{$A_2=A_3=0$ for all solutions} \\ \hline
\end{tabular}
}
\end{table}

\begin{table}[tbp]
\caption{Braided Fibonacci oscillator of the second type}
\label{table3}\centering
{\scriptsize
\begin{tabular}{|c||c|c|}
\hline
\multicolumn{3}{|c|}{$aa^{*}-Q_1a^{*}a=q^{N}$,\ $aq^{N}=qq^Na$,\ $%
q^{N}a^{*}=qa^{*}q^{N}$} \\ \hline
& sol1 & sol2 \\ \hline
$A_2=$ & $\displaystyle \frac{q-Q_1}{q}$ & $\displaystyle \frac{Q_1-q}{Q_1}$
\\ \hline
$A_3=$ & $\displaystyle \frac{Q_1(Q_1-q)}{q}$ & $q-Q_1$ \\ \hline
$k_1=$ & $-\displaystyle \frac{Q_1}{q}$ & $-\displaystyle \frac{q}{Q_1}$ \\
\hline
$b_2=$ & $0$ & $\displaystyle \frac{Q_1^2-q^2}{qQ_1}$ \\ \hline
$b_3=$ & $q$ & $q$ \\ \hline
$z=$ & $\displaystyle \frac{q}{Q_1}$ & $\displaystyle \frac{Q_1}{q}$ \\
\hline
$c_2=$ & $q^{-1}$ & $q^{-1}$ \\ \hline
$c_3=$ & $\displaystyle \frac{q^2-Q_1^2}{qQ_1}$ & $0$ \\ \hline
$d_1=$ & $\displaystyle \frac{1}{Q_1}$ & $\displaystyle \frac{Q_1}{q^2}$ \\
\hline
$f_1=$ & $\displaystyle \frac{q^2}{Q_1}$ & $Q_1$ \\ \hline
$g_1=$ & $\displaystyle \frac{q^2}{Q_1^2}$ & $\displaystyle \frac{Q_1^2}{q^2}
$ \\ \hline
\multicolumn{3}{|c|}{$b_1=c_1=d_2=f_2=g_2=g_3=0$} \\
\multicolumn{3}{|c|}{for both solutions} \\ \hline
\end{tabular}
}
\end{table}
\begin{table}[tbp]
\caption{Other braided oscillator solutions}
\label{table5}\centering
{\scriptsize
\begin{tabular}{|c||c|c|c|c|c|c|c|c|}
\hline
\multicolumn{6}{|c|}{$aa^{*}-Q_1a^{*}a=Q_2q^{2N}+Q_3q^N$,\ $aq^{N}=qq^Na$,\ $%
q^{N}a^{*}=qa^{*}q^{N}$} &  &  &  \\ \hline
& sol1 & sol2 & sol3 & sol4 & sol5 &  &  &  \\ \hline
$Q_1=$ & $q^2$ & $-q^2$ & $-q^2$ & $q^2$ & $-q$ &  &  &  \\ \hline
$Q_2=$ & $1$ & $1$ & $1$ & $1$ & $0$ &  &  &  \\ \hline
$Q_3=$ & $0$ & $0$ & $0$ & $0$ & $1$ &  &  &  \\ \hline
$k_1=$ & $-1$ & $-1$ & $-1$ & $-1$ & $-1$ &  &  &  \\ \hline
$b_1=$ & $0$ & $0$ & $2-q^2c_1$ & $2+q^2c_1$ & $0$ &  &  &  \\ \hline
$b_3=$ & $q^{2}$ & $-q^{2}$ & $-q^{2}$ & $q^{2}$ & $-q$ &  &  &  \\ \hline
$z=$ & $1$ & $-1$ & $-1$ & $1$ & $-1$ &  &  &  \\ \hline
$c_1=$ & $0$ & $0$ & $c_1$ & $c_1$ & $0$ &  &  &  \\ \hline
$c_2=$ & $q^{-2}$ & $-q^{-2}$ & $-q^{-2}$ & $q^{-2}$ & $-q^{-1}$ &  &  &  \\
\hline
$d_1=$ & $q^{-1}$ & $q^{-1}$ & $q^{-1}$ & $q^{-1}$ & $q^{-1}$ &  &  &  \\
\hline
$f_1=$ & $q$ & $q$ & $q$ & $q$ & $q$ &  &  &  \\ \hline
$g_1=$ & $-1$ & $-1$ & $1$ & $1$ & $1$ &  &  &  \\ \hline
\multicolumn{6}{|c|}{$A_2=A_3=b_2=c_3=d_2=f_2=g_2=g_3=0$} &  &  &  \\
\multicolumn{6}{|c|}{for all solutions} &  &  &  \\ \hline
\end{tabular}
}
\end{table}

\begin{table}[tbp]
\caption{Three dimensional braided space}
\label{table4}\centering
{\scriptsize
\begin{tabular}{|c||c|c|c|c|c|c|c|}
\hline
\multicolumn{8}{|c|}{$aa^{*}=qa^{*}a$,\ $aq^{N}=qq^Na$,\ $%
q^{N}a^{*}=qa^{*}q^{N}$} \\ \hline
& sol1 & sol2 & sol3 & sol4 & sol5 & sol6 & sol7 \\ \hline
$A_2=$ & $A_2$ & $A_2$ & $0$ & $0$ & $0$ & $0$ & $0$ \\ \hline
$A_3=$ & $-qA_2$ & $qA_2$ & $0$ & $0$ & $0$ & $0$ & $0$ \\ \hline
$k_1=$ & $-1$ & $1$ & $-1$ & $-1$ & $-1$ & $-1$ & $-1$ \\ \hline
$b_2=$ & $0$ & $0$ & $0$ & $0$ & $-1+g_1$ & $-1+qc_2$ & $0$ \\ \hline
$b_3=$ & $q$ & $q$ & $qg_1$ & $q$ & $q$ & $q$ & $qz$ \\ \hline
$z=$ & $1$ & $1$ & $g_1$ & $z$ & $g_1$ & $qc_2$ & $z$ \\ \hline
$c_2=$ & $q^{-1}$ & $q^{-1}$ & $q^{-1}$ & $q^{-1}$ & $q^{-1}g_1$ & $c_2$ & $%
q^{-1}$ \\ \hline
$c_3=$ & $0$ & $0$ & $g_1-1$ & $0$ & $0$ & $0$ & $z-1$ \\ \hline
$d_1=$ & $q^{-1}$ & $q^{-1}$ & $q^{-1}$ & $q^{-1}$ & $q^{-1}g_1$ & $q^{-1}$
& $q^{-1}$ \\ \hline
$d_2=$ & $0$ & $0$ & $g_1-1$ & $0$ & $0$ & $0$ & $0$ \\ \hline
$f_1=$ & $q$ & $q$ & $qg_1$ & $q$ & $q$ & $q$ & $q$ \\ \hline
$f_2=$ & $0$ & $0$ & $0$ & $0$ & $-1+g_1$ & $0$ & $0$ \\ \hline
$g_1=$ & $1$ & $1$ & $g_1$ & $g_1$ & $g_1$ & $g_1$ & $g_1$ \\ \hline
\multicolumn{8}{|c|}{$b_1=c_1=g_2=g_3=0$ for all solutions} \\ \hline
\end{tabular}
}
\end{table}

\end{document}